\documentclass[MSNbibl,number,citesort,dvips]{arxbj}
\usepackage{mathrsfs}
\usepackage{graphicx}

\aid{0}
\volume{20}
\issue{4}
\pubyear{2014}
\firstpage{1999}
\lastpage{2019}
\doi{10.3150/13-BEJ548} 

\makeatletter
\mathchardef\varTheta="0102
\newcommand{\eqref}[1]{(\ref{#1})}
\renewcommand{\epsilon}{\varepsilon}
\newcommand{\binom}[2]{{#1}\choose {#2}}
\newcommand{\red}[1]{#1}
\renewcommand{\d}{\mathrm{d}}
\newcommand{\X}{\mathcal{X}}
\newcommand{\Y}{\mathcal{Y}} 
\newcommand{\M}{\mathcal{M}} 
\newcommand{\Th}{\varTheta} 
\newcommand{\nuo}{\nu} 
\newcommand{\meas}{\red{\xi}}
\newcommand{\Z}{\mathbb{Z}_+}
\newcommand{\Du}{D} 
\renewcommand{\r}{r} 
\renewcommand{\O}{\mathcal{O}}
\newcommand{\N}{\mathbb{N}}
\newcommand{\R}{\mathbb{R}}
\newcommand{\E}{\mathbb{E}}
\newcommand{\F}{\mathcal{F}}
\renewcommand{\L}{\mathcal{L}} 
\newcommand{\Gx}{\mathcal{A}} 
\newcommand{\Gm}{A} 
\newcommand{\Fmix}{\bar{\F}_f}
\newcommand{\D}{\mathscr{D}}
\newcommand{\bs}[1]{\boldsymbol{#1}}
\renewcommand{\aa}{\bs\alpha}
\newcommand{\ee}{\mathbf{e}}
\newcommand{\mm}{\mathbf{m}}
\newcommand{\nn}{\mathbf{n}}
\newcommand{\ii}{\mathbf{i}}
\newcommand{\oo}{\mathbf{0}}
\renewcommand{\Pr}{\mathbb{P}}
\newcommand{\hw}{\widehat{w}}

\newtheorem{proposition}{Proposition}[section]
\newtheorem{lemma}[proposition]{Lemma}
\newtheorem{Lemma}{Lemma}[section]
\makeatother

\begin{document}
\begin{frontmatter}

\title{Optimal filtering and the dual process}
\runtitle{Filtering and the dual}

\begin{aug}
\author[1]{\inits{O.}\fnms{Omiros} \snm{Papaspiliopoulos}\corref{}\thanksref{1}\ead[label=e1]{omiros.papaspiliopoulos@upf.edu}} \and
\author[2]{\inits{M.}\fnms{Matteo} \snm{Ruggiero}\thanksref{2}\ead[label=e2]{matteo.ruggiero@unito.it}}
\address[1]{ICREA \& Department of Economics and Business, Universitat
Pompeu Fabra,
Ram\'on Trias Fargas 25-27, 08005, Barcelona, Spain. \printead{e1}}
\address[2]{Collegio Carlo Alberto
\& Department of Economics and Statistics, University of Torino,\\  C.so Unione Sovietica 218/bis, 10134,
Torino, Italy. \printead{e2}}
\end{aug}

\received{\smonth{5} \syear{2013}}

%
\begin{abstract}
We link optimal filtering for hidden Markov models to the notion of
duality for Markov processes. We show that when the signal is dual
to a process that has two components, one deterministic and
one a pure death process, and with respect to functions that define
changes of measure conjugate to the emission density, the
filtering distributions evolve in the family of finite mixtures of such
measures and the filter can be computed at a cost that is polynomial
in the number of observations.
Special cases
of our framework include the Kalman filter, and computable filters
for the Cox--Ingersoll--Ross process and the one-dimensional Wright--Fisher
process, which have been investigated before.
The dual we obtain for the Cox--Ingersoll--Ross process appears to be new
in the literature.
\end{abstract}

%
\begin{keyword}
\kwd{Bayesian conjugacy}
\kwd{Cox--Ingersoll--Ross process}
\kwd{finite mixture models}
\kwd{hidden Markov model}
\kwd{Kalman filter}
\end{keyword}

\end{frontmatter}

\section{Introduction}
\label{sec:intro}

\begin{figure}

\includegraphics{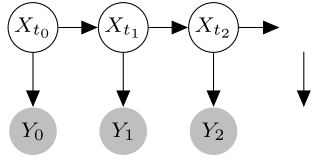}

\caption{Hidden Markov model represented as a graphical model.}\label{hmm}
\end{figure}

A hidden Markov model (HMM) for a sequence of observations $\{Y_n,n\geq
0\}$, where $Y_n \in\Y$, is a discrete-time stochastic process with
dynamics depicted
in Figure~\ref{hmm}. It is defined in terms of a hidden Markov
chain, the so-called
\emph{signal}, which in this paper will be taken to be the
discrete-time sampling of
a time-homogeneous continuous-time Markov process $X_t$, with
state-space $\X$,
transition kernel $P_t(x, \mathrm{d} x')$, and initial distribution $\nuo(\d
x)$.
The observations relate to the signal by means of conditional
distributions, assumed to be given by the kernel $F(x, \d y)$. We will
assume that
%
\begin{equation}
\label{eq:emission} F(x,\d y) = f_x(y) \mu(\d y)
\end{equation}
for some measure $\mu(\d y)$, in which case the corresponding densities
are known
as the \emph{observation} or \emph{emission}
densities.
The optimal filtering
problem is the derivation of the
conditional distributions $\L(X_{t_n} | Y_0,\ldots,Y_n)$ of the
unobserved signal given the observations
collected up to time $t_{n}$, henceforth denoted $\nu_n(\d x)$. These filtering
distributions are the backbone of all statistical
estimation problems in this framework, such as the prediction of
future observations, the derivation of smoothing distributions (i.e.,
the conditional distribution of $X_{t_n}$ given past and future
observations) and the calculation of the likelihood function, that is,
the marginal density of the observations when the emission
distributions are dominated.
See
\cite{Cea05} for details and applications.

Throughout the paper, we will assume that the signal is stationary and
reversible with respect to a probability measure $\pi$.
\red{Section~\ref{sec:disc} shows how to extend our result to
non-stationary signals.}
It is also
appealing, from a modeling point of view,
to assume that the signal evolves in continuous time, since there is a
rich family of such models with a prespecified stationary measure
$\pi$. In addition, this assumption will give us a powerful tool to
study optimal
filtering by using the generator of the process, as we show
in Section~\ref{sec:dual}. In the examples of
Section~\ref{sec:examples}, the state space $\X$ of the
signal will either be a
subset of $\R$ or the
$(K-1)$-dimensional simplex $\Delta_K$.

Mathematically, optimal
filtering is the solution of the recursion
\[
\nu_0 = \phi_{Y_0} (\nuo), \qquad \nu_{n} =
\phi_{Y_{n}}\bigl(\psi_{ t_n -
t_{n-1}}(\nu_n)\bigr),\qquad n>0,
\]
which involves the following two operators acting on probability
measures $\meas$:
%
\begin{equation}
\label{update-predict} %
\begin{array} {@{}l@{\qquad}l@{}} \mbox{update:} &
\phi_{y}(\meas) (\d x) =\displaystyle\frac{f_x(y) \meas(\d x)}{p_{\meas}(y)},\qquad
p_{\meas}(y) = \int_{\X} f_x(y) \meas(\d
x),
\\[9pt]
\mbox{prediction:}& \psi_t(\meas) \bigl(\d x'\bigr)=
\meas P_t\bigl(\d x'\bigr)=\displaystyle\int
_{\X}\meas(\d x)P_t\bigl(x,\d x'
\bigr). \end{array} %
\end{equation}
The ``update'' is the application of Bayes theorem, and the
``prediction'' gives the distribution of the next step of the Markov
chain initiated from $\meas$.
These operators have
the following property when applied to
\emph{finite mixtures of distributions}:
%
\begin{equation}
\label{linearities} \phi_{y} \Biggl(\sum_{i=1}^n
w_i \meas_i \Biggr) (\d x) =\sum
_{i=1}^n {w_i p_{\meas_i}(y) \over\sum_j w_j
p_{\meas_j}(y)}
\phi_y(\meas_i), \qquad \psi_t \Biggl(\sum
_{i=1}^n w_i
\meas_i \Biggr) (\d x) =\sum_{i=1}^n
w_i \psi_t(\meas_i).
\end{equation}
This implies that when $\X$ is a finite set, there is a
simple algorithm for the sequential computation of the filtering
probabilities. To see this, note that we can think of a
distribution $\nu$ on a finite set $\X$, specified in terms of
probabilities $\alpha_x, x \in\X$, as a finite mixture of point
masses, $\nu= \sum_x \alpha_x \delta_x$; it is easy
to compute $\phi_y(\delta_x),\psi_t(\delta_x)$ and then use the above
result to
obtain the probabilities associated with the distributions
$\phi_y(\nu)$ and $\psi_t(\nu)$. This yields a popular algorithm
for inference in HMMs, commonly known as the
Baum--Welch filter, whose complexity is easily seen to be $\O(n |\X|^2)$,
where $|\X|$ is the cardinality of $\X$.

Outside the finite state-space case, the iteration of
these two operators typically leads to
analytically intractable distributions. However, there are notable
exceptions to this rule. The
classic example is the linear Gaussian state-space model, for
which the filtering distributions are Gaussian with mean and
covariance that can be iteratively computed using the so-called
Kalman filter, at cost that grows linearly with $n$. Recent work by
Genon-Catalot and collaborators uncovered that there exist
interesting non-Gaussian models for which the filtering distributions are
\emph{finite mixtures} of\vadjust{\goodbreak} parametric distributions. See \cite
{GK04,CG06,CG09}, where the authors show how to compute the
corresponding parameters sequentially in
these models. We revisit their findings in Section~\ref{sec:examples}. However, the number of mixture components increases
with $n$ in a way such that the cost
of computing the filters grows polynomially with $n$ (see Section~\ref{sec:dual} for details). Borrowing and
adapting the terminology from \cite{CG06}, we will refer to
filters with such computational cost as
\emph{computable}, whereas filters whose cost grows linearly with $n$
as \emph{finite-dimensional}.

The work by Genon-Catalot and collaborators raises four important questions,
which we address in this paper:
are there more models which admit computable filters;
do they share some basic structure; is there a general methodology to
identify such models
and to obtain the algorithm which computes the sequence of parameters;
what is the computational complexity of such schemes and how can we
obtain faster approximate filtering algorithms? We show that the answer
to all
these questions relates to an important probabilistic object: the \emph
{dual process}.
Duality methods have a long history in Probability, dating back
to the work of P. L\'evy \cite{L48} (see \cite{JK13} for a recent
review). These have been widely applied to the study of interacting
particle systems \cite{L05} and proven to be a powerful method which
provides alternative, and often simpler, tools for investigating the
sample path properties of the process at hand. For example, the
existence of a dual for a certain Markov process (and for a
sufficiently large class of functions) implies that the associated
martingale problem is well defined, hence that the process is
unique; see Section~4.4 of \cite{EK86}. See also \cite{D93} and
\cite{E00} for applications of duality to population genetics.

In this paper, we illustrate that dual processes play a central role in
optimal filtering and to a great extent can be used to settle the
four questions posed above. We also uncover their potential as
auxiliary variables in Monte
Carlo schemes for stochastic processes (and, hence, as a
variance reduction scheme). In our framework, the dual will in general
be given by two components: a \emph{deterministic process}, driven by
an ordinary differential equation, and a (multidimensional)
\emph{death process} with countable state-space. We show how to derive
an explicit, recursive filtering scheme once the dual is identified,
and apply this methodology to three cases of fundamental interest.
In doing so, we identify what, to the best of our
knowledge, is a new \textit{gamma-type}
duality.

The rest of the paper is organized as follows. In Section~\ref{sec:dual}, we link optimal filtering to a specific type of
duality, we show how to identify the dual in terms of the
\emph{generator} of $X_t$, and study the complexity of the resulting
filtering algorithm.
Section~\ref{sec:examples} analyzes three interesting models for which
the dual process is derived:
the Cox--Ingersoll--Ross model, the Ornstein--Uhlenbeck process and the
$K$-dimensional Wright--Fisher diffusion. These models are reversible
with respect to
the gamma, Gaussian and Dirichlet distribution,
respectively, and for the Gaussian case the computable
filter reduces to the Kalman filter. Section~\ref{sec:disc} discusses
certain aspects of the methodology, including the extension to
infinite-dimensional signals modeled as Fleming--Viot processes.

\section{Methodology: Filtering the dual process}
\label{sec:dual}

\subsection{Linking optimal filtering to duality}

Before presenting the main results, we introduce three fundamental
assumptions which provide the general framework under which the
results are derived.
First, we will assume that $X$ is
reversible with respect to a probability measure $\pi$:
\begin{itemize}
\item[A1 (Reversibility):] $\pi(\d x) P_t(x, \d x') = \pi(\d x')
P_t(x', \d x)$.
\end{itemize}
\red{Section~\ref{sec:disc} discusses how this assumption can be
relaxed to accommodate non-stationary signals.}
In order to state the second assumption, we need to introduce a certain
amount of
notation. Define, for $K \in \Z=\N\cup\{0\}$, the space of multi-indices
%
\begin{equation}
\label{eq:M} \M=\Z^K= \bigl\{ \mm=(m_1,
\ldots,m_K)\dvtx m_j \in\Z, j=1, \ldots, K \bigr\}.
\end{equation}
We will use the symbol $\oo$ to denote the vector of zeros, $\ee_j$ for
the vector in $\M$
whose only non-zero element is found at the $j$th coordinate and equals
$1$, and let $|\mm|=\sum_im_i$. Furthermore, we will use the product
order on $\M$, according to which for
$\mm, \nn\in\M$, $\mm\leq\nn$ if an only if $m_j \leq n_j$ for all
$j$. Then, for $\ii\leq\mm$, $\mm-\ii$ is the vector with $j$th
element $m_j-i_j$. Additionally, if $\Lambda\subset\M$, define
%
\begin{equation}
\label{eq:gl} G(\Lambda) = \{\red{\nn\in\M}\dvtx \nn\leq\mm, \mm\in\Lambda\}.
\end{equation}
The notation for $\M$ does not reflect its dependence on
the dimension $K$, but we will reserve boldface for elements of $\M$
when $K>1$ (or unspecified), whereas normal typeface will be
used for elements of $\Z$. Finally, the following notations will be
used to denote conditional expectations
%
\[
(P_t f) (x) = \E^x\bigl[f(X_t)
\bigr] = \E\bigl[f(X_{t})\vert X_{0}=x\bigr] = \int
_\X f\bigl(x'\bigr) P_t\bigl(x,
\d x'\bigr).
\]
The first denotes the action on $f$ of the
semigroup operator associated to the
transition kernel, where with some
abuse of notation the same symbol is used both for the semigroup and
the kernel.

The second assumption is concerned with models where
$\pi(\d x)$ is \emph{conjugate} to the emission density $f_x(y)$:
\begin{itemize}
\item[A2 (Conjugacy):]
For $\Th\subseteq\R^l, l \in\Z$, let $h\dvtx \X\times\M\times\Th\to
\R_+$ be such that $\sup_x
h(x,\mm,\allowbreak  \theta) < \infty$ for all $\mm\in\M,\theta\in\Th$, and
$h(x,\oo,\tilde{\theta}) = 1$ for some $\tilde{\theta}
\in\Th$.
Then $\F=\{h(x,\mm,\theta) \*\pi(\d x), \mm\in\M, \theta\in\Th\}$ is
assumed to be a
family of probability measures such that there exist functions
$t\dvtx \Y\times\M\to\M$ and $T\dvtx  \Y\times\Th\to\Th$ with $\mm\to
t(y,\mm)$ increasing and
such that
\[
\phi_y\bigl(h(x,\mm,\theta) \pi(\d x)\bigr) = h\bigl(x,t(y,
\mm),T(y,\theta)\bigr) \pi(\d x).
\]
\end{itemize}
Hence here with conjugacy, we intend the fact that the family $\F$ of
measures, which includes $\pi$, is
closed under the update operation. The
assumption that $h$ is bounded in $x$ will be discussed after the
statement of Assumption A4.

For $p_\nu(y)$ as in
\eqref{update-predict}, it is easy to check that in the context of
A2, we have
%
\begin{equation}
\label{eq:norm-const} p_{h(x,\mm,\theta)\pi(\d x)}(y) =: c(\mm,\theta,y) = {f_x(y)
h(x,\mm,\theta) \over h(x,t(y,m),T(y,\theta))},
\end{equation}
which, despite its appearance, does not depend on $x$.

Note that our definitions of $\M$ and $\Th$ allow the possibility that $K=0$
or $l=0$, in which case $h$ in A2 is function only of the
variables with non-zero dimension, whereas the case $K=l=0$ is not of interest
here.
In the setting of Assumption A2 and for the trivial Markov dynamics
$X_t\equiv X_0$, with $X_0 \sim\pi$, the filtering problem collapses
to conjugate
Bayesian inference for the unknown parameter $x$ of the sampling density
$f_x(y)$. See Section~5.2 and Appendix A.2 of \cite{BSbook}
for an exposition of conjugate Bayesian inference and stylized
conjugate Bayesian models, and Section~\ref{sec:examples} in this
paper for examples within our
framework.

The third main assumption for our results concerns the existence of a
certain type of dual process for the signal.
\begin{itemize}
\item[A3 (Duality):]
We assume that $\r\dvtx \Th\to\Th$ is such that the differential equation
%
\begin{equation}
\label{eq:ode} \d\Theta_t / \d t = \r(\Theta_t),\qquad
\Theta_0 = \theta_0,
\end{equation}
has a unique solution for all $\theta_0$.
Let
$\lambda\dvtx  \Z \to\R_+$ be an increasing function, $\rho\dvtx \Th\to
\R_+$ be a continuous function, and consider a two-component Markov
process $(M_t,\Theta_t)$ with state-space $\M\times\Th$, where $\Theta
_t$ evolves autonomously
according to \eqref{eq:ode}, and when at
$(M_t,\Theta_t)=(\mm,\theta)$, the process jumps down to state
$(\mm-\ee_j,\theta)$ with instantaneous rate
%
\begin{equation}
\label{eq:dual-rate} \lambda\bigl(|\mm|\bigr) \rho(\theta) m_j.
\end{equation}
We assume $(M_t,\Theta_t)$ is \emph{dual} to $X_t$ with respect to
the family of functions $h$ defined in A2, in the sense that
%
\begin{equation}
\label{eq:means} \E^x\bigl[h(X_t,\mm,\theta)\bigr] =
\E^{(\mm,\theta)}\bigl[h(x,M_t,\Theta_t)\bigr]\qquad
\forall x \in\X, \mm\in\M, \theta\in\Th, t\geq0.
\end{equation}
%
When $K=0$ or $l=0$ in A2, the dual process is just $\Theta_t$ or
$M_t$, respectively, and we adopt the convention that
\[
\rho(\theta)\equiv1 \quad\mbox{whenever}\quad l=0.
\]
\end{itemize}

Note that $M_t$ can only jump to ``smaller''
states according to the partial order on $\M$,
and that \eqref{eq:dual-rate} implies that 0 is an absorbing state for
each coordinate $j$ of $M_{t}$, so that the vector of zeros is a global
absorbing state.

As mentioned in Section~\ref{sec:intro}, the notion of duality for
Markov processes with respect to a given function is well known. See,
for example, Section II.4 in \cite{L05}.
Among the most common type of duality relations we mention \emph{moment
duality}, that is duality with respect to functions of type
$h(x,y)=x^{y}$, and \emph{Laplace duality}, that is with respect to
functions of type $h(x,y)=\mathrm{e}^{-axy}$. See, for example, \cite{JK13}.
In our framework, the duality functions are Radon--Nikodym derivatives
between measures
that are conjugate to the emission density, and this setup is
perfectly tailored to optimal filtering.
Furthermore, A3 specifies that we are interested in
dual processes which can be decomposed into two parts: one
purely \emph{deterministic} and the other given by a $K$-dimensional
\emph{pure death process}, whose death rates are subordinated by the
deterministic process. The transition probabilities of the death
process, conditional on the initial state $\Theta_{0}=\theta$, will be
denoted by
%
\begin{equation}
\label{eq:trans-dual} p_{\mm,\nn}(t;\theta) = \Pr[M_t = \nn|
M_0=\mm,\Theta_0=\theta], \qquad \nn, \mm\in\M, \nn\leq
\mm.
\end{equation}
It is worth mentioning that the requirements on the structure of the
dual processes prescribed by Assumption A3, with particular reference
to the
intensity \eqref{eq:dual-rate}, are\vadjust{\goodbreak} justified by the three main
reasons. The first is that, as shown in Section~\ref{sec:examples},
they define a framework general enough to identify
duals of processes of interest, the incorporation of a deterministic
component being necessary in this respect. The second reason is that
the transition probabilities \eqref{eq:trans-dual} are analytically
available, as provided by the following result, whose proof can be
found in the \hyperref[sec: proofs]{Appendix}.

\begin{proposition}\label{prop:trans probab}
Let $(M_{t},\Theta_t)$ be as in \textup{A}3, with $(M_{0},\Theta_0)=(\mm,\theta
)\in\M\times\Th$, and let
$\lambda_{|\mm|}=|\mm|\lambda(|\mm|)$. Then the transition
probabilities for $M_{t}$ are
$p_{\mm,\mm}(t;\theta)=\exp\{-\lambda_{|\mm|}\int_{0}^{t}\rho(\Theta
_{s})\,\d s\}$
and, for any $\oo\le\ii\le\mm$,
\[
p_{\mm,\mm-\ii}(t;\theta)= \Biggl(\prod
_{h=0}^{|\ii|-1}\lambda_{|\mm
|-h} \Biggr)
C_{|\mm|,|\mm|-|\ii|}(t) p\bigl(i_{1},\ldots,i_{K}; \mm,|\ii|\bigr),
\]
where
\[
C_{|\mm|,|\mm|-|\ii|}(t)= (-1)^{|\ii|}\sum
_{k=0}^{|\ii|}\frac{\mathrm{e}^{-\lambda_{|\mm|-k}\int
_{0}^{t}\rho(\Theta_{s})\,\d s}}{\prod_{0\le h\le|\ii|,h\ne k}(\lambda
_{|\mm|-k}-\lambda_{|\mm|-h})}
\]
and $p(i_{1},\ldots,i_{K}; \mm,|\ii|)$ is the multivariate
hypergeometric probability mass function with parameters $(\mm,|\ii|)$
evaluated at $(i_{1},\ldots,i_{K})$.
\end{proposition}

This result can be interpreted as follows. The probability that a
one-dimensional death process with inhomogeneous rates
$\lambda_{|\mm|}\rho(\Theta_{s})$ decreases from $|\mm|$ to $|\mm|-|\ii
|$ in the\vspace*{-1pt} interval
$[0,t]$ is
$(\prod_{h=0}^{|\ii|-1}\lambda_{|\mm|-h})C_{|\mm|,|\mm|-|\ii|}(t)$,
where the second factor is related to the convolution of the waiting
times in an inhomogeneous Poisson process (see Section~19.10 in \cite
{JK94}, and \cite{SB99}).
For a $K$-dimensional death process, the same quantity is the probability
associated to all paths leading from level $|\mm|$ to level
$|\mm|-|\ii|$.
Given such event, the probability of the subset of paths leading
exactly from $\mm$ to $\mm-\ii$ is then given by the multivariate
hypergeometric probability $p(i_{1},\ldots,i_{K}; \mm,|\ii|)$ (an
expression of this probability can be found in the \hyperref[sec: proofs]{Appendix}).

Note that the special case of Proposition~\ref{prop:trans probab}
yielded by $K=1$ and $\rho(\Theta_{s})\equiv1$ relates to the result
obtained in Proposition~4.5 in
\cite{CG09}. Note also that when $\rho(\Theta_{s})\equiv1$ and
$\lambda_{m}=m(\theta+m-1)/2$, $C_{|\mm|,|\mm|-|\ii|}(t)$ is the
transition probability of the block-counting process of Kingman's
coalescent with mutation, see \cite{T84} and
\cite{G06} for details on such process.

The third motivation behind the type of duality required by A3 is
that if it holds, the prediction operator maps measures as in A2 into
finite mixtures.
%
\begin{proposition}
\label{prop:prox-mix}
Let $\psi_{t}$ be as in (\ref{update-predict}) and assume \textup{A}1--\textup{A}2--\textup{A}3
hold. Then
%
\begin{equation}
\label{eq:propagation} \psi_t\bigl(h(x,\mm,\theta) \pi(\d x)\bigr) = \sum
_{\oo\leq\ii\leq\mm} p_{\mm,\mm-\ii}(t;\theta) h(x,\mm-\ii,
\Theta_t) \pi(\d x)
\end{equation}
with $p_{\mm,\mm-\ii}(t;\theta)$ as in Proposition~\ref{prop:trans
probab} and where $\Theta_t$ is the value in $t$ of the process
in \eqref{eq:ode}
started from $\Theta_0=\theta$.
\end{proposition}

\begin{pf}
From \eqref{update-predict}, we have
\begin{eqnarray*}
\psi_t\bigl(h(x,\mm,\theta) \pi(\d x)\bigr) &=& \int
_\X h(x,\mm,\theta) \pi(\d x) P_t\bigl(x,\d
x'\bigr) = \int_\X h(x,\mm,\theta) \pi\bigl(
\d x'\bigr) P_t\bigl(x',\d x\bigr)
\\
&=& \pi\bigl(\d x'\bigr) \E^{x'}\bigl[h(X_t,
\mm,\theta)\bigr] = \pi\bigl(\d x'\bigr) \E^{(\mm,\theta)}\bigl[h
\bigl(x',M_t,\Theta_t\bigr)\bigr]
\\
&=& \sum_{\nn\leq\mm} p_{\mm,\nn}(t;\theta) h
\bigl(x',\nn,\Theta_t\bigr) \pi\bigl(\d x'
\bigr),
\end{eqnarray*}
\red{where the second equality follows from A1, the fourth from A3, and
the last from \eqref{eq:trans-dual}.}
\end{pf}

The above result states that reversibility and the existence of the
required duality
jointly guarantee that the prediction operator can be computed with a
finite effort. The reduction of the operator to a sum is due to the
fact that $X_t$ is dual to a Markov process with discrete state-space, but
it is precisely the fact that $M_t$ is a pure death process that makes
the number of terms in the sum being finite. The next result shows
that computable filtering is available in the framework we have
outlined.
%
\begin{proposition}
\label{prop:comp}
Consider the family of finite mixtures
%
\begin{equation}
\label{eq:fm} \Fmix= \biggl\{\sum_{\mm\in\Lambda}
w_{\mm} h(x,\mm,\theta) \pi(\d x) \dvtx \Lambda\subset\M, |\Lambda| <
\infty, w_\mm\geq0, \sum_{\mm\in\Lambda}
w_\mm=1 \biggr\}.
\end{equation}
Then, under Assumptions \textup{A}1--\textup{A}2--\textup{A}3, $\Fmix$ is closed under the
application of the prediction and update operators \eqref
{update-predict}, and specifically
\[
\phi_y \biggl( \sum_{\mm\in\Lambda}
w_{\mm} h(x,\mm,\theta) \pi(\d x) \biggr) = \sum
_{\nn\in t(y,\Lambda)} \hw_\nn h\bigl(x,\nn,T(y,\theta)\bigr) \pi(
\d x)
\]
with
%
\begin{eqnarray}
\label{eq:update-mix} t(y,\Lambda) &:=& \bigl\{\nn: \nn=t(y,\mm), \mm\in\Lambda\bigr\}
\nonumber
\\[-8pt]
\\[-8pt]
\hw_\nn &\propto& w_{\mm} c(\mm,\theta,y) \qquad\mbox{for }
\nn= t(y,\mm), \sum_{\nn
\in t(y,\Lambda)} \hw_\nn=1,
\nonumber
\end{eqnarray}
and
%
\begin{equation}
\label{eq:predict-mix} \psi_t \biggl( \sum_{\mm\in\Lambda}
w_{\mm} h(x,\mm,\theta) \pi(\d x) \biggr) = \sum
_{\nn\in G(\Lambda)} \biggl( \sum_{{\mm\in\Lambda, \mm\geq\nn}}
w_\mm p_{\mm,\nn}(t;\theta) \biggr) h(x,\nn,
\theta_t) \pi(\d x).
\end{equation}
\end{proposition}

The above proposition shows that under Assumption
A1 to A3, and provided the starting state belongs to the family $\Fmix$
of finite mixtures with
components as in A2, then the filtering
distributions evolve within $\Fmix$. Furthermore, the explicit
reweighing of the
mixture components is provided, thus allowing to concretely implement
the recursive filtering scheme.
Note also that this result generalizes Theorem~2.1 in \cite{CG06},
which states a similar
result for $K=1$ under the Assumption A2 and the
result in Proposition~\ref{prop:prox-mix}.
The proof of Proposition~\ref{prop:comp} follows from
\eqref{linearities}, A2 and Proposition~\ref{prop:prox-mix} by direct
computation, and is thus omitted. Later in this section, we will
derive filtering algorithms based on this result. However, we first
address in the next subsection the most important aspect of the
approach described in
this section, which is how
to find a dual process that satisfies A2.

\subsection{Local duality as a sufficient condition}

It is typically easier to identify
a process that satisfies the duality relation \eqref{eq:means} for infinitesimal
$t$. Formally, this requires studying the \emph{generator} of $X_t$,
which we will denote by
$\Gx$. This is a linear
operator, with domain denoted $\D(\Gx)$,
linked to the semigroup operator via the Kolmogorov
backward equation
%
\[
{\partial\over\partial t} P_t f(x) = (\Gx
P_t f) (x),\qquad f \in\D(\Gx ),
\]
where on the left hand side $P_t h(x)$ is differentiated in $t$
for given $x$, whereas on
the right hand side, $\Gx$ acts on $P_t h(x)$ as a function of $x$
for given $t$. See, for example, Proposition~1.1.5 in \cite{EK86}.

Suppose now $X_t$ is a diffusion process which solves an
SDE on $\R^d$ of the form
%
\[
\d X_t =b(X_t) \,\d t +
\sigma(X_t) \,\d B_t.
\]
In this case, $\Gx$ is the second-order differential operator given by
%
\begin{equation}
\label{eq:gen-dif} (\Gx f) (x) = \sum_{i=1}^{d}b_{i}(x)
\frac{\partial f(x)}{\partial
x_{i}} + \frac{1}{2} \sum_{i,j=1}^{d}a_{i,j}(x)
\frac{\partial^{2}
f(x)}{\partial x_{i}\,\partial x_{j}},\qquad f\in\D(\Gx),
\end{equation}
for an appropriate domain $\D(\Gx)$ and where $a_{i,j}(x) := (\sigma(x)
\sigma(x)^\mathrm{T})_{i,j}$.

Let now $\Gm$ denote the generator of the dual process defined
in A3, which can be easily checked to be
%
\begin{eqnarray}\label{eq:gen-dual}
 (\Gm g) (\mm,\theta) &=& \lambda\bigl(|\mm|\bigr) \rho(\theta)
\sum_{i=1}^K m_i \red{
\bigl[g(\mm-\ee_i,\theta)-g(\mm,\theta)\bigr]}
\nonumber
\\[-8pt]
\\[-8pt]
&&{}+ \sum
_{i=1}^{l}\r _{i}(\theta)
\frac{\partial g(\mm,\theta)}{\partial\theta}, \qquad g\in\D (\Gm),
\nonumber
\end{eqnarray}
with $\r$ as in \eqref{eq:ode}.
The main idea is then to identify the dual process from the generator,
instead of
the semigroup operator.
\begin{itemize}
\item[A4 (Local duality):] The function $h(x,\mm,\theta)$ defined in
A2 is such that $h(x,\mm,\theta)$, as a function of $x$ belongs to
$\D(\Gx)$ for all $(\mm,\theta) \in\M\times\Th$, as a function of
$(\mm,\theta)$ belongs to $\D(\Gm)$ for all $x \in\X$, and
%
\begin{equation}
\label{eq:gen-eq} \bigl(\Gx h(\cdot,\mm,\theta)\bigr) (x) = \bigl(\Gm h(x,\cdot,
\cdot)\bigr) (\mm,\theta) \qquad \forall x \in\X, \mm\in\M, \theta\in\Th.
\end{equation}
\end{itemize}
It will typically be the case that in order to have $h(\cdot,\mm,\theta
) \in\D(\Gx)$, one needs
$h$ to be bounded in $x$, requirement included in Assumption A2.
Proposition~1.2 of \cite{JK13} shows that A4, together with the further
assumption
\[
\E^x\bigl[h(X_t,\mm,\theta)\bigr]
\in\D(\Gm),\qquad \E^{(\mm,\theta)}\bigl[h(x,M_t,
\Theta_t)\bigr]
\in\D(\Gx),
\]
implies \eqref{eq:means}, whose argument can be sketched as follows.
From (\ref{eq:gen-eq}), we can write
\[
(\beta I-\Gx) h = (\beta I-\Gm) h,\qquad \beta\in\R,
\]
where $I$ denotes the identity operator.
Since $\Gx$ and $\Gm$ generate strongly continuous contraction
semigroups, say on $L_{1}$ and $L_{2}$, their ranges are dense in
$L_{1}$ and $L_{2}$, respectively. Moreover, the resolvents
$\mathcal{R}_\beta=(\beta I-\Gx)^{-1}$, $R_\beta=(\beta I-\Gm)^{-1}$ are
one-to-one for all $\beta>0$, so the previous implies
\[
\mathcal{R}_\beta h=R_\beta h,\qquad h\in
L_{1}\cap L_{2}, \beta>0.
\]
Since the resolvent of an operator is the
Laplace transform of the associated semigroup, and because of the
uniqueness of Laplace transforms, the previous expression
in turn implies \eqref{eq:means}.

The approach sketched above for identifying the dual process by means
of the local condition \eqref{eq:gen-eq} will be implemented in
Section~\ref{sec:examples}, where we will identify the duals for some
interesting relevant models.

\subsection{The filtering algorithm}\label{sec:alg}

\begin{figure}

\includegraphics{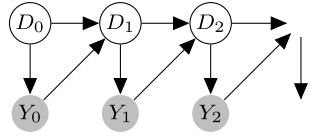}

\caption{The partially observed Markov process dual to the hidden
Markov model in
Figure~\protect\ref{hmm}, where $\Du_i=(M_i,\Theta_i)$.}\label{hmm2}
\end{figure}

Typically, the initial distribution of the signal process belongs to $\F
$, and most often equals the invariant measure $\pi$. Thus, without
loss of generality and in order to simplify the exposition below, we
make the following additional
assumption.
\begin{itemize}
\item[A5 (Initialization):] The initial distribution of the signal is
$\nu= h(x,\mm_0,\theta_0) \pi(\d x) \in\F$, for some
$\mm_0 \in
\M$, $\theta_0 \in\Th$.
\end{itemize}
Proposition~\ref{prop:comp} provides a probabilistic
interpretation of the weights involved in the finite mixtures in terms
of the transition probabilities of the dual death process $M_{t}$. This
interpretation can be elaborated further, in order to facilitate the
development of filtering algorithms. With a little abuse of notation,
denote by $\{\Du_n=(M_n,\Theta_n), n\geq0\}$ a discrete-time
process with
state-space $\M\times\Th$ constructed as follows. Consider a partially observed
Markov process, where the signal is now $\Du_n$
and the conditional independence structure, given in Figure~\ref{hmm2} graphically, is as follows.
Let $D_0
=(M_0,\Theta_0)=(\mm_0,\theta_0)$ be the initial state of the chain,
with $(\mm_0,\theta_0)$ defined in A5.
Then $\L(Y_n | D_n=(\mm,\theta)) = c(\mm,\theta,y) \mu(\d
y)$, with $c(\mm,\theta,y) $ as in \eqref{eq:norm-const} and $\mu$ in
\eqref{eq:emission}, and for $n \geq1$, $\L(D_n | Y_{n-1}=y,
D_{n-1}=(\mm,\theta))$ is the law of
$(M_{t_n-t_{n-1}},\Theta_{t_n-t_{n-1}})$ in A3 started from
$(t(y,\mm),T(y,\theta))$ at time 0.
Then, the connection between duality and optimal filtering can be
expressed as
%
\begin{equation}
\label{eq:dual-filter} \L(X_{t_n} | Y_0,\ldots,Y_n) =
\int h\bigl(x,t(Y_n,M_n),T(Y_n,
\Theta_n)\bigr) \pi (\d x) \,\d \L(D_n | Y_0,
\ldots,Y_{n-1}).
\end{equation}
Thus, filtering $X_{t_i}$ in the original model in Figure~\ref{hmm}
can be achieved by filtering $D_i$ in the dual model in Figure~\ref{hmm2}. Since $\Theta_n$
evolves deterministically, optimal filtering for $X$ reduces to filtering
$M_{n}$, which has finite support with probabilities that can be
computed recursively using an algorithm similar to the Baum--Welch
filter, as
we now describe. $\L(M_0,\Theta_0)$ has
support on the single point $\{(\mm_0,\theta_0)\}$;
if $\L(M_n,\Theta_n |
Y_0,\ldots,Y_{n-1})$ has support on $\Lambda_n \times\{\theta_n\}$
for $\Lambda_n \subset\M$ and $\theta_n \in\Th$, and assigns probability
$w_\mm$ to state $(\mm,\theta_n)$, then $\L(M_{n+1},\Theta_{n+1} |
Y_0,\ldots,Y_n)$
has support on $\Lambda_{n+1} \times\{\theta_{n+1}\}$, where
$\Lambda_{n+1} = G(t(Y_n,\Lambda_n))$, for $G$ and $t(y,\cdot)$ defined
in
\eqref{eq:gl} and \eqref{eq:update-mix}, respectively, $\theta_{n+1}$
the solution of \eqref{eq:ode} at time $t_{n+1}-t_{n}$ started from
$\Theta_0=\theta_n$, and the
probability associated to state $(\nn,\theta_{n+1}) \in\Lambda_{n+1}
\times
\{\theta_{n+1}\}$ is
%
\begin{equation}
\label{eq:dual-prop} \Pr[M_{n+1} =\nn, \Theta_{n+1}=
\theta_{n+1} | Y_0,\ldots,Y_n] = \mathop{\sum
_{\mm\in
\Lambda_{n}}}_ {t(Y_n,\mm)\geq\nn} w_\mm
p_{t(Y_n,\mm),\nn}(t_{n+1}-t_n; \theta_n).
\end{equation}
Therefore, the optimal filtering reduces to the sequential computation
of the parameters $\theta_n$, the supports $\Lambda_n$ and the
probabilities on each support point in $\Lambda_{n}$, for
$n=0,1,\ldots$\,. 

%

The computation of the probabilities \eqref{eq:dual-prop} for all $(\nn
,\theta_{n+1}) \in\Lambda_{n+1} \times
\{\theta_{n+1}\}$ can be done at a cost that is at most of order
$|\Lambda_{n+1}|^2$. Therefore, the overall cost of computing the
filters up to the $n$th observation is bounded from above by $\sum_{i=0}^n
|\Lambda_i|^2$. If $|\Lambda_i|$ were constant with $i$, we would
recover the complexity of the Baum--Welch filter, discussed in Section~\ref{sec:intro}. However, $|\Lambda_i|$ increases with $i$, as a
result of the successive operation of $G$ and $t(y,\cdot)$ defined
in \eqref{eq:gl} and \eqref{eq:update-mix}, respectively. Clearly, it is
hard to make further analysis on the computational complexity without
some information on
$t(y,\cdot)$. Here we will assume that $t(y,\mm)=\mm+N(y)$, where
$N\dvtx \Y\to\M$, a structure that is found in all the examples we
study in this paper. We then have the following key result. The proof
of the lemma is omitted.
%
\begin{lemma}
\label{prop:g}
For any $\Lambda\subset\M$ and $\mm\in\M$, we define $\Lambda+ \mm=
\{\nn+\mm; \nn\in\Lambda\}$. Then
\[
G\bigl(G(\Lambda)+\mm\bigr) = G(\Lambda+ \mm).
\]
\end{lemma}

\begin{proposition}\label{prop:complex} Under the assumption that $t(y,\mm)=\mm+N(y)$, where
$N\dvtx \Y\to\M$, we have that
\[
|\Lambda_n| = G \Biggl(\mm_0+\sum
_{i=1}^n N(Y_i) \Biggr) \leq \biggl( 1+
{d_n \over K} \biggr)^K,
\]
where $d_n = |\mm_0+\sum_{i=1}^n N(Y_i)|$.
\end{proposition}
\begin{pf}
The equality follows by successive application of Lemma~\ref{prop:g}. For the inequality, notice that $\nn:= \mm_0+\sum_{i=1}^n
N(Y_i) \in\M$, with $|\nn|=d_n$. Then, by \eqref{eq:gl}, $|G(\nn)|
=\prod_{i=1}^K(n_i+1)$. Then, apply Jensen's inequality to $\log
|G(\nn)|$ to obtain the result.
\end{pf}

When the observations follow a stationary process, $d_n$ will be of
order $n$. Therefore, the complexity of carrying out the computations
involved in the filtering recursions up to iteration $n$, will be
\red{$\mathcal{O}(n^{2K})$}, where the constant depends on $K$ but not
$n$. We return to the issue of complexity in Section~\ref{sec:disc}.


\section{The dual of some stochastic processes}
\label{sec:examples}

Following the local duality approach outlined in the previous section,
here we identify the dual processes for the Cox--Ingersoll--Ross
model, diffusion processes with linear coefficients and $K$-dimensional
Wright--Fisher
diffusions. \red{In accordance with the rest of the article, we focus
on stationary parametrizations of the processes and discuss the
non-stationary case in Section~\ref{sec:disc}.}

Recall that $d,K,l$ denote the state space dimension for $X_{t}$, $M_{t}$
and $\Theta_{t}$, respectively.

\subsection{CIR processes}\label{ex:lsde}

The so-called Cox--Ingersoll--Ross (CIR) model is a non-negative one-dimensional
diffusion, that solves the SDE
%
\[
\d X_t =\bigl(\delta\sigma^{2}-2\gamma
X_{t}\bigr) \,\d t + 2\sigma\sqrt{X_{t}} \,\d
B_t.
\]
This name is due to \cite{CIR85} who
introduced the model in mathematical finance, although this model had been
studied long before in the literature, see, for example, the population
growth model in Section~13.C of \cite{KT81} and the process described in Section~5 of
\cite{F51}. From a broader perspective, the CIR model can also be seen
as a special case
of a continuous-state branching process with immigration \cite{KW71}.

The generator of the CIR process is
%
\begin{equation}
\label{eq:CIR} \Gx=\bigl(\delta\sigma^{2}-2\gamma x\bigr)
\frac{\d}{\d x}+ 2\sigma^{2}x\frac{\d^{2}}{\d x^{2}},\qquad \delta,\gamma,
\sigma>0,
\end{equation}
with domain defined as follows. With the above
parametrization, and using Feller's terminology, the boundary point
$+\infty$ is natural for all choices of parameters, while 0 is regular if
$\delta<2$ and entrance if $\delta\ge2$. Define
\[
\D_{0}(\Gx)= \bigl\{f\in C_{0}\bigl([0,\infty
)\bigr)\cap C^{2}\bigl((0,\infty)\bigr)\dvtx \Gx f\in C_{0}
\bigl([0,\infty)\bigr) \bigr\},
\]
where $C_{0}([0,\infty))$ is the space of continuous functions
vanishing at infinity,
and
\[
\D(\Gx)= \cases{ %
f\in\D_{0}(\Gx),
& \quad\mbox{if } $\delta\ge2$,\vspace*{2pt}
\cr
f\in\D_{0}(\Gx)\dvtx
\displaystyle\lim_{x\rightarrow0}x^{\delta/2}f'(x)=0, &
\quad\mbox{if } $0<\delta<2$. 
}
\]
Then $\{(f,\Gx f)\dvtx f\in\D(\Gx)\}$ generates a Feller
semigroup on $C_{0}([0,\infty))$.
\red{Such choice of the domain for the case $0<\delta<2$ guarantees
that the boundary 0 is instantaneously reflecting.}
See Theorems 8.1.2 and 8.2.1 in \cite{EK86}.
In this case,
the CIR process is reversible with respect to the gamma distribution
\[
\pi\equiv\operatorname{Ga}\bigl(\delta/2, \gamma/\sigma^{2}\bigr).
\]
Previous results on duality for the CIR model include a Laplace
duality, that is with respect to a function of type $h(x,y)=\mathrm{e}^{-axy}$.
See, for example, \cite{HW07}.
Here however we identify a new, gamma-type duality relation, which
has as special cases a moment and a Laplace duality.
Now, let $d=1$, $K=1$ and $l=1$, and
define, for $\theta>0$, the function
\[
h(x,m,\theta)= \frac{\Gamma(\delta/2)}{\Gamma(\delta/2+m)} \red{\biggl({\gamma\over
\sigma^{2}}
\biggr)^{-\delta/2}}\theta^{\delta/2+m}x^{m} \exp\bigl\{-\bigl(
\theta-\gamma/\sigma^{2}\bigr)x\bigr\}.
\]
This function can be identified as the Radon--Nikodym derivative of a
$\operatorname{Ga}(\delta/2+m,\theta)$ distribution with
respect to $\pi$. The family of gamma distributions that arises
by varying $m \in\Z$ and $\theta>0$, defines a subset of the family of gamma
distributions that is conjugate to emission
densities that as a function of $x$ are proportional to
\[
x^n \mathrm{e}^{-\lambda x}, \qquad n \in\Z, \lambda>0,
\]
in which case $t$ and $T$ in A2 coincide with
\[
t(y,m) = n+m, \qquad T(y,\theta) = \theta+\lambda.
\]
Such type of emission density arises, for example, for observations
$Y_n=n$ distributed as
Poisson with intensity $\lambda X_{t_n}$, giving rise to a
dynamic version of the Poisson-gamma conjugate Bayesian model.

On the other hand, $h(\cdot,m,\theta)$ belongs to the domain of $\Gx$ only
when $\theta\geq\gamma/\sigma^2$, in which case $h\in
C^{2}_{0}([0,\infty))$. In
order to be able to use local duality as in A4 we will assume that the
family is defined as
\[
\F=\bigl\{h(x,m,\theta) \pi(\d x), m \in\Z, \theta\geq \gamma/
\sigma^2\bigr\}
\]
but we will return to the case $\theta<\gamma/\sigma^2$ at the end of this
subsection.
Then a simple computation yields
\begin{eqnarray*}
\Gx h(\cdot,m,\theta) (x) &=& 2m\sigma^2\theta h(x,m-1,
\theta)+\sigma^2 (\delta+2m) \bigl(\theta-\gamma /\sigma^2
\bigr)h(x,m+1,\theta)
\\
&&{} -\sigma^2\bigl[2m\theta+(\delta+2m) \bigl(\theta-\gamma/
\sigma^2\bigr)\bigr]h(x,m,\theta).
\end{eqnarray*}
Motivated by this structure, and with view to achieving the local
duality in \eqref{eq:gen-eq} we consider a two-component process
$(M_{t},\Theta_{t})$ with generator $\Gm$ as in \eqref{eq:gen-dual}, where
\[
\red{\lambda(m)=2\sigma^{2}}, \qquad \r(\theta) = 2
\sigma^2 \theta\bigl(\gamma/\sigma^2-\theta\bigr), \qquad
\rho(\theta)=\theta.
\]
It is then easy to check that local duality holds, namely
\[
\Gx h(\cdot,m,\theta) (x)=\Gm h(x,\cdot,\cdot) (m,\theta).
\]
Additionally, the conditions that are required to derive
\eqref{eq:means} from this local duality are satisfied.
In this example, the solution of the dynamical system \eqref{eq:ode}
for $\Theta_0=\theta$
is given by
\[
\Theta_t = {\gamma\over\sigma^2} {\theta \mathrm{e}^{2\gamma t} \over
\theta \mathrm{e}^{2\gamma t} + \gamma/\sigma^2 - \theta},
\]
which in conjunction with Proposition~\ref{prop:trans probab} implies
that the
transition probabilities for the death process simplify to binomial
probabilities
\[
p_{m,m-i}(t;\theta)=\operatorname{Bin} \biggl(m-i ; m,
{\gamma\over
\sigma^2} \bigl( \theta \mathrm{e}^{2\gamma t} + \gamma/
\sigma^2 - \theta\bigr)^{-1} \biggr).
\]
Therefore, we have all the ingredients necessary to implement the
filtering algorithm. Finally, note that if $\theta_0\geq\gamma/\sigma
^2$, then
$\theta_n\geq\gamma/\sigma^2$ for all $n$.

Notice that the result on the transition probabilities above, together
with Proposition~\ref{prop:prox-mix}, implies the following interesting
property of the CIR process:
%
\begin{eqnarray}
\label{eq:mixture-CIR} &&\psi_t\bigl( \operatorname{Ga} (m+\delta/2, \theta)
\bigr)
\nonumber
\\[-8pt]
\\[-8pt]
&&\quad =\sum_{k=0}^m \operatorname{Bin}
\biggl(k ; m, {\gamma\over
\sigma^2} \bigl( \theta \mathrm{e}^{2\gamma t} +
\gamma/\sigma^2 - \theta\bigr)^{-1} \biggr)
\operatorname{Ga} \biggl(k+\delta/2, {\gamma\over\sigma^2} {\theta \mathrm{e}^{2\gamma t} \over
\theta \mathrm{e}^{2\gamma t} + \gamma/\sigma^2 - \theta}
\biggr).
\nonumber
\end{eqnarray}
This result has been obtained before, using a completely
different approach; the case
$\delta=1$ can be shown directly by
elementary calculations using a change of variables and binomial
expansion of the left-hand-side; the general case was proved in
\cite{CG06}, see Proposition~3.4 and the associated Lemma~3.1, after
some rather heavy calculations. The result in \eqref{eq:mixture-CIR}
leads to a computable
filter, as we showed in Proposition~\ref{prop:comp}, which is
precisely the result also obtained in \cite{CG06} for the CIR
process. It is neat that using duality
and the generic result in Proposition~\ref{prop:prox-mix}, this result
can be obtained in a straightforward manner. The proof in \cite{CG06}
is based on the following known series expansion of the CIR transition
kernel, see expression (80) in \cite{CG06} and page 334 of \cite{KT81},
which can be re-expressed as a Poisson mixture of gamma distributions
as follows:
%
\begin{equation}
\label{eq:series} P_t\bigl(x,\d x'\bigr) = \sum
_{k \geq0} \operatorname{Poisson} \biggl(k ; {\gamma
\over\sigma^2}
{1 \over \mathrm{e}^{2\gamma t} -1} x \biggr) \operatorname{Ga} \biggl(k+\delta/2,
{\gamma
\over\sigma^2} {\mathrm{e}^{2\gamma t} \over \mathrm{e}^{2\gamma t} -1} \biggr).
\end{equation}
It is interesting that instead of deriving \eqref{eq:mixture-CIR} from
\eqref{eq:series}, which in any case is laborious, one can prove the
\red{former using duality and then obtain \eqref{eq:series} by taking
$\theta=(m+\delta/2)/x$ and letting $m\to\infty$ in
\eqref{eq:mixture-CIR}.}

In view of the arguments of Section~\ref{sec:dual}, it follows that,
for $\theta< \gamma/\sigma^2$, $h\notin\D(\Gx)$, hence duality in the
sense of
A3 cannot be established using local duality. However, in
view of the result \eqref{eq:mixture-CIR} that has already been
obtained in \cite{CG06}, it is obvious that duality still holds in
this case. This also shows the limitation of the functional analytic
method for establishing duality: it is a very powerful when all formal
requirements are met, but
there will be examples, like this one, where \eqref{eq:means} would
have to be established by alternative arguments. Nevertheless, a formal
calculation using the generator reveals the dual even when $\theta<
\gamma/\sigma^2$.

\subsection{Linear diffusion processes}\label{ex:lsde}

We consider the scalar Ornstein--Uhlenbeck process that
solves an SDE of the form
%
\[
\d X_t =-{\sigma^2 \over\alpha}(X_t -
\gamma) \,\d t + \sqrt{2} \sigma\,\d B_t,
\]
which is reversible with respect to the Gaussian distribution,
\[
\pi(\d x) \equiv\operatorname{Normal}(\gamma,\alpha).
\]
The generator is given by
\[
\Gx= \bigl(\sigma^2 \gamma/\alpha- \sigma^2 x/\alpha
\bigr) {\d\over\d x} + \sigma^2 {\d^2 \over\d x^2}
\]
with domain $C_0^2((-\infty,\infty))$. In this model, we have $d=1,
K=0,l=2$, where $\theta=(\mu,\tau) \in\R\times\R_{+}$, and
\[
h(x,\mu,\tau) = \biggl({\alpha\over\tau} \biggr)^{1/2} \exp
\biggl\{ -{(x-\mu)^2 \over2 \tau} + {(x-\gamma)^2 \over2 \alpha} \biggr\},
\]
which can be easily recognised as the Radon--Nikodym derivative between
a $\operatorname{Normal}(\mu,\tau)$ and~$\pi$. The measures $h(x,\mu,\tau)
\pi(\d x)$ are conjugate to emission
densities that as a function of $x$ are proportional to
\[
\exp \biggl\{-{1 \over2 \lambda} (x-c)^2 \biggr\}, \qquad
\lambda>0, c \in\R,
\]
with $T(y,\theta)=((\lambda\mu+ \tau c), \lambda\tau)/
(\lambda+\tau)$. Such density arises, for example, with data $Y_n=c$
that is Gaussian with mean $X_{t_n}$ and variance $\lambda$.
As with the CIR process, we
have the technical problem that this function belongs to $\D(\Gx)$
only for $\tau<\alpha$, hence we will restrict to this case and define
$\Th=\{(\mu,\tau)\dvtx  \mu\in\R, 0<\tau<\alpha\}$.
A direct calculation gives that
\[
\Gx h(\cdot,\mu,\tau) (x) = {\sigma^2 \over\alpha} (\gamma-\mu)
{\partial\over\partial\mu} h(x,\mu,\tau)+ 2 \sigma^2(1-\tau/\alpha)
{\partial\over\partial\tau} h(x,\mu,\tau).
\]
This suggests that the dual is purely deterministic and described in terms
of the ODEs:
\[
\d\mu_t / \d t = {\sigma^2 \over\alpha} (\gamma-
\mu_t) \,\d t, \qquad \d\tau_t / \d t = 2
\sigma^2(1-\tau_t/\alpha) \,\d t.
\]
Duality with respect to this deterministic process implies that the
filter evolves within the Gaussian family and the computational cost
is linear in $n$, that is, we are dealing with a finite-dimensional filter.

Of course, all this is known: the ODEs above are the well-known
equations for the first two moments of linear SDEs, and the filter is
the Kalman filter. Thus, within the assumptions we have made in this
article, the finite-dimensional filter corresponds to the special case
where the dual is purely deterministic. We considered $d=1$ for
simplicity, but the results carry over to multi-dimensional stationary
linear SDEs. The same discussion as for the CIR applies here regarding
the restrictions posed by needing that $h \in\D(\Gx)$. We return to
this issue in Section~\ref{sec:disc}.

\subsection{Wright--Fisher diffusions} \label{ex:wf}

Wright--Fisher (WF) processes are $K$-dimensional diffusions with paths
confined in the $(K-1)$-dimensional simplex
%
\begin{equation}
\label{eq:K-simplex} \Delta_{K}= \Biggl\{x\in[0,1]^{K}\dvtx \sum
_{i=1}^{K}x_{i}=1 \Biggr\}.
\end{equation}
These processes approximate,
among others, large-population
discrete Wright--Fisher
reproductive models with non-overlapping generations, and describe the
time-evolution of the species abundancies when the individuals in the
underlying population are subject to random genetic drift and,
possibly, mutation, selection and recombination. See, for example,
Chapter~5 in \cite{D93} and Chapter~10 in \cite{EK86}.
Here we are interested in the case without
selection nor recombination, and with \emph{parent-independent
mutation}. \red{That is, we consider a WF diffusion with generator
%
\begin{equation}
\label{eq:WF-generator} \Gx= \frac{1}{2}\sum_{i=1}^{K}\bigl(
\alpha_{i}-|\aa| x_{j}\bigr)\frac{\partial
}{\partial x_{i}} +\frac{1}{2}
\sum_{i,j=1}^{K}x_{i}(
\delta_{ij}-x_{j})\frac{\partial
^{2}}{\partial x_{i}\,\partial x_{j}},
\end{equation}
where $\delta_{ij}$ denotes the Kronecker delta, $\aa=(\alpha_{1},\ldots
,\alpha_{K})\in\R_{+}^{K}$ and $|\aa|=\sum_{i=1}^{K}\alpha_{i}$. The
domain of the operator $\Gx$ is taken to be
$C^{2}(\Delta_{K})$, and the closure of $\Gx$ generates a strongly
continuous contractive semigroup on $C(\Delta_{K})$. See \cite{EK81}
for details.}
Note that this is a hypoelliptic diffusion, that is, the square of the
diffusion matrix is not full rank, as a result of the constraint
$\sum_i x_i=1$. Even though we could work with an elliptic diffusion
for the $K-1$
variables, it is the formulation above that is desirable for
identifying the dual, as we will show.

Such diffusion is reversible with respect to the Dirichlet
distribution
%
\begin{equation}
\label{eq:Dirichlet} \pi(\d x_{1},\ldots,\d x_{K}) =
\frac{\Gamma(|\aa|)}{\prod_{j=1}^{K}\Gamma(\alpha_{j})}x_{1}^{\alpha
_{1}-1}\cdots x_{K}^{\alpha_{K}-1}
\,\d x_1 \,\cdots\,\d x_K,\qquad x\in
\Delta_{K}.
\end{equation}

In this model, we have $d=K\ge2$ and $l=0$, therefore
there is no deterministic component in the dual process. We denote
\[
x^{\mm}=x_{1}^{m_{1}}\cdots x_{K}^{m_{K}},
\qquad x\in\Delta_{K}, \mm\in\M,
\]
and define
%
\begin{equation}
\label{eq:C3} h(x,\mm)=\frac{\Gamma(|\aa|+|\mm|)}{\Gamma(|\aa|)}\prod_{j=1}^{K}
\frac
{\Gamma(\alpha_{j})}{\Gamma(\alpha_{j}+m_{j})}x^{\mm},
\end{equation}
whence clearly $h(\cdot,\mm)\in\D(\Gx)$. This can be identified with
the Radon--Nikodym derivative between a Dirichlet distribution with
parameters $(\alpha_1+m_1,\ldots,\alpha_K+m_K)$ and $\pi$, and it is
conjugate to emission densities that as a function of $x$ are
proportional to
\[
x_{1}^{n_{1}}\cdots x_{K}^{n_{K}},
\qquad n_{i} \in\Z, i=1,\ldots,K,
\]
in which case $t$ in A2 coincides with $t(y,m) = n+m$.
Such type of emission density arises, for example, for observations
$Y_n=(n_{1},\ldots,n_{K})$ distributed as
Multinomial with parameters $X_{t_n}=(X_{t_{n},1},\ldots,X_{t_{n},K})$,
giving rise to a
dynamic version of the Dirichlet-Multinomial conjugate Bayesian model.

Then we have
\begin{eqnarray*}
\Gx h(x,\mm) &=& \sum
_{i=1}^{K} \biggl(\frac{\alpha_{i}m_{i}}{2}+{
\binom{m_{i}} {2}} \biggr)\frac{\Gamma(|\aa|+|\mm|)}{\Gamma(|\aa|)}\prod
_{j=1}^{K}\frac{\Gamma
(\alpha_{j})}{\Gamma(\alpha_{j}+m_{j})}x^{\mm-\ee_{i}}
\\
&&{} -\sum_{i=1}^{K} \biggl(
\frac{|\aa| m_{i}}{2}+{\binom{m_{i}} {2}}+\frac
{1}{2}m_{i}
\sum_{j\ne i}m_{j} \biggr)\frac{\Gamma(|\aa|+|\mm|)}{\Gamma
(|\aa|)}
\prod_{j=1}^{K}\frac{\Gamma(\alpha_{j})}{\Gamma(\alpha
_{j}+m_{j})}x^{\mm}
\\
&=& \frac{|\aa|+|\mm|-1}{2}\sum_{i=1}^{K}m_{i}h(x,
\mm-\ee_{i})-\frac{|\mm
|(|\aa|+|\mm|-1)}{2}h(x,\mm). 
\end{eqnarray*}
This suggests considering a one-component dual process, with $M_{t}$ a
Markov jump process with generator $A$ obtained by letting
\[
\lambda\bigl(|\mm|\bigr)=\bigl(|\aa|+|\mm|-1\bigr)/2,\qquad \rho(\theta)\equiv1,
\]
in (\ref{eq:gen-dual}).
Since $h(x,\cdot)\in\D(\Gm)$, it is then easy to check that the local
duality condition
\[
\Gx h(\cdot,\mm) (x)= \Gm h(x,\cdot) (\mm)
\]
holds.
Hence, the WF diffusion with parent-independent mutation $X_{t}$ and
the death process $M_{t}$ on $\Z^{K}$, which jumps from $\mm$ to $\mm
-\ee_{j}$ at rate $m_{j}(|\aa|+|\mm|-1)/2$, are dual with respect to
the above $h$ in the sense of A3. The transition probabilities of
$M_{t}$ are as in Proposition~\ref{prop:trans probab}.

Filtering the WF model when $K=2$ on the basis of binomial
data was studied in \cite{CG09}. One can appreciate the strength of
the approach we introduce here, since it is straightforward to obtain
the filtering recursion using the dual and Proposition~\ref{prop:trans
probab} for any $K$. It has to be noted that, in our opinion, one of
the reasons why the results are harder to obtain using the approach in
\cite{CG09}, is because they decide to work with the elliptic
WF model, which is a scalar
diffusion since $K=2$.
Working with the elliptic model hides the structure of
duality, which is immediately apparent in the hypoelliptic model.

\red{The death process we obtain in this section can be seen as a
special case of the process used in \cite{BEG00} for deriving an
infinite mixture expansion for the transition kernel of the WF
diffusion with selection}.

An extension of WF diffusions to the case of infinitely-many types is
given by Fleming--Viot processes. These are measure-valued diffusions
whose finite-dimensional projections onto partitions of the type space
coincide with WF processes.
A duality relation holds between the Fleming--Viot process and a
function-valued process related to Kingman's coalescent. See, for
example, \cite{EK93}.
However, by applying to such dual process the same finite-dimensional
projection that yields the WF process, one does not obtain the dual
derived here, since binning the process into finitely-many sets hides
some important information about the events at the level of particles.

\section{Discussion}
\label{sec:disc}

We have demonstrated that computable filtering follows from duality,
in the sense described in Assumptions A2 and A3 in
Section~\ref{sec:dual}. A sufficient condition to establish duality is the
local duality described in Assumption A4, which is based on the
properties of the generator of the signal process and its relation to
the semigroup operator via the Kolmogorov backward equation. The use
of this functional analytic machinery places some constraints on the
duality function in A2, such as for example that as a function of $x$
it has to vanish at infinity. Therefore, even when duality holds in
the sense of A3 for functions that do not satisfy such constraints,
the local duality cannot be used to prove this. On the other hand, the
local duality can still be used formally to identify the dual. Both in
the CIR process when $\theta<\gamma/\sigma^2$ and in the OU process
when $\tau<\alpha$ (see Section~\ref{sec:examples} for details) the
formal application of the generator identifies the dual correctly.

We have assumed reversibility with respect to a probability measure
$\pi$.
\red{In fact, our methodology relies on the existence of such
reversible measure but does not require that it be a probability
measure. An inspection of Proposition~\ref{prop:prox-mix} reveals that
it is still valid in this more general case, provided $h(x,\mm,\theta
)\pi(\d x)$ is a probability measure. Therefore, $h$ is the
Radon--Nikodym derivative between the measures in $\F$ and the
reversible measure, in this more general framework which also covers
non-stationary signals.}

Another topic of investigation is the connection of the duality, as
used in this paper, and results about the spectral representation of
the transition kernel of the signal, for example, the type of
expression in \eqref{eq:series}. There are classic results about such
expressions, see, for example, Chapter~13 of \cite{KT81}, and their existence
seems to be related to computable filtering, see Section~6.4 of
\cite{CG09}, but the connection is not well understood.

Our results in Section~\ref{sec:dual} show that for observations
generated by a stationary process the computational cost associated
with the identification of the filtering distributions grows
polynomially with the number of observations, unless $K=0$ in which
case the growth is linear. However, it might be the case that most of
the components in the mixture representations have negligible
weight. Previous simulation studies, see, for example, Table~1 in
\cite{GK04}, show that after a few iterations the filter might
concentrate all its mass in two or three components. We believe that
the connection to the dual process might be very
helpful in studying the effective number of components. However, there
are subtleties in this line of research. Note that when $X_t$ is
ergodic, and $t_i-t_{i-1}$ is large relative to its mixing time,
practically all mass of the filtering distribution will be
concentrated on a single component, the ``root''
$(\oo,\tilde{\theta})$ (see A2) that corresponds to the invariant
measure $\pi$. Therefore, the
time evolution of the number of states with non-negligible filtering
probabilities (say above a given $\epsilon\approx0$) will depend on
the number of observations per unit of time in the $X$ process. This
aspect deserves careful study.


We have obtained explicit filters also for a class of Fleming--Viot diffusions with
parent independent mutation and for a class of measure-valued branching processes
with immigration, under parametrizations which make them reversible with respect
to the Dirichlet and the gamma process respectively. A peculiarity of this
framework, entailed by the fact that the signal is measure-valued and the
observations are
random draws from the signal, is the lack of a common dominating measure for the
emission distributions, hence the lack of likelihood, which makes the nature of the
problem truly non-parametric. The techniques for obtaining an optimal filter thus
necessarily differ from those illustrated in this paper and will be reported elsewhere.

\begin{appendix}\label{sec: proofs}
\section*{Appendix: Proof of Proposition \texorpdfstring{\protect\ref{prop:trans probab}}{2.1}}

Before stating the result, we recall a useful lemma, whose proof can be
found in \cite{CG09}.

\begin{Lemma}\label{lemma CG09}
\[
\sum_{j=0}^{l}
\frac{(-1)^{j}}{(\lambda_{n}-\lambda_{n-1-j})\prod_{0\le
h\le l,h\ne j}|\lambda_{n-1-j}-\lambda_{n-1-h}|} =\frac{1}{\prod_{1\le h\le l+1}(\lambda_{n}-\lambda_{n-h})}.
\]
\end{Lemma}

\begin{pf*}{Proof of Proposition~\ref{prop:trans probab}}
Consider first
the one-dimensional case, that is, $\mm=m$, denote for brevity
$\theta[s,t]=\int_{s}^{t}\rho(\theta_{u})\,\d u$, and define for $i\ge1$
\begin{eqnarray*}
I_{1,\ldots,i}&=& \int_{0}^{t}
\cdots\int_{t_{i-1}}^{t}\mathrm{e}^{-\lambda
_{m}\theta[0,t_{1}]} \prod
_{k=1}^{i-1}\theta_{t_{k}}
\mathrm{e}^{-\lambda_{m-k}\theta
[t_{k},t_{k+1}]}\,\d t_{k}\,\theta_{t_{i}}
\mathrm{e}^{-\lambda_{m-i}\theta
[t_{i},t]}\,\d t_{i},
\\
I_{1,\ldots,j-1,j,\ldots,i} &=& \int_{0}^{t}\cdots\int
_{t_{i-1}}^{t}\mathrm{e}^{-\lambda_{m}\theta[0,t_{1}]} \prod
_{k=1,k\ne j}^{i-1}\theta_{t_{k}}
\mathrm{e}^{-\lambda_{m-k}\theta
[t_{k},t_{k+1}]}\,\d t_{k} \,\theta_{t_{i}}
\mathrm{e}^{-\lambda_{m-i}\theta[t_{i},t]}\,\d t_{i},
\end{eqnarray*}
where $t_{j}:=t_{j+1}$ in $I_{1,\ldots,j-1,j+1,\ldots,i}$. It can be
easily seen that
\setcounter{equation}{0}
\begin{equation}
\label{eq:I-i} \red{I_{i}=\frac{\mathrm{e}^{-\lambda_{m-i}\theta[0,t]}-\mathrm{e}^{-\lambda_{m}\theta
[0,t]}}{\lambda_{m}-\lambda_{m-i}}}.
\end{equation}
Then we have
%
\begin{equation}
\label{eq:rescale-p-i} \Biggl(\prod_{h=0}^{i-1}
\lambda_{m-h} \Biggr)^{-1}p_{m,m-i}(t)=I_{1,\ldots,i},
\end{equation}
where $p_{m,m-i}(t)$ is the transition probability associated \red{to
the one-dimensional death process.}
By integrating twice, we obtain
\begin{eqnarray*}
I_{1,\ldots,i} &=& \frac{(-1)(I_{1,\ldots,i-1}-I_{1,\ldots,i-2,i})}{\lambda
_{m-(i-1)}-\lambda_{m-i}}
\\
&=& \frac{(-1)^{2}}{\lambda_{m-(i-1)}-\lambda_{m-i}} \biggl[ \frac{(I_{1,\ldots,i-2}-I_{1,\ldots,i-3,i-1})}{\lambda_{m-(i-2)}-\lambda
_{m-(i-1)}} -\frac{(I_{1,\ldots,i-2}-I_{1,\ldots,i-3,i})}{\lambda_{m-(i-2)}-\lambda
_{m-i}} \biggr].
\end{eqnarray*}
The iteration of the successive integrations can be represented as a
binary tree with root $(i,0):=I_{1,\ldots,i}$, whose node
$(i-j,i-k):=I_{1,\ldots,i-j,i-k}$ branches into $(i-j,0):=I_{1,\ldots
,i-j}$ and $((i-j-1)^{+},i-k)=I_{1,\ldots,i-j-1,i-k}$, with both
branches weighed $1/(\lambda_{m-(i-j)}-\lambda_{m-(i-k)})$, determined
by the parent node's indices. The leaves correspond to nodes where the
left coordinate touches zero if the right coordinate is already zero,
or where the left crosses zero if the right coordinate is positive. The
term associated to the leaf $(0,i-k)$ will be $(-1)^{i}\mathrm{e}^{-\lambda
_{m-(i-k)}\theta[0,t]}$ weighed by some appropriate coefficient.
%
%
The level before the leaves can be seen as the sequence
\[
\underbrace{\underbrace{\underbrace{\underbrace{I_{1}I_{2}}_{{2^{1}}}
\underbrace{I_{1}I_{3}}}_{{2^{2}}}
\underbrace{I_{1}I_{2}I_{1}I_{4}}}_{{2^{3}}}
\underbrace{I_{1}I_{2}I_{1}I_{3}I_{1}I_{2}I_{1}I_{5}}}_{{2^{4}}}
\ldots,
\]
where every sequence of $2^{i}$ terms is repeated with the last index
augmented by one,
and each $I_{i}$ produces the leaves $\mathrm{e}^{-\lambda_{m}}$ and $\mathrm{e}^{-\lambda
_{m-i}}$. Hence given $i$, there are $2^{i-2}$ terms $I_{1}$, $2^{i-3}$
terms $I_{2}$, \ldots, $2^{1}$ terms $I_{i-2}$, $2^{0}$ terms $I_{i-1}$
and $I_{i}$.
Note also that $I_{1}$ has
$2^{0}$ paths in common with $I_{i}$,
$2^{0}$ paths in common with $I_{i-1}$,
$2^{1}$ paths in common with $I_{i-2}$,
\ldots,
$2^{i-3}$ paths in common with $I_{2}$.

The correct coefficient for $I_{k}$ is computed by collecting some
constants related to the paths that have the same last coefficient and
simplifying. In particular, given $i$, the paths to be grouped for
$I_{k}$ are those whose constants for indices greater than $k$
change, since according to the rule above, when $k$ is the rightmost
index in $I_{1,\ldots,k}$, there is only one path down to $I_{k}$.
Hence, given $i$, term $I_{k}$ has coefficient
\begin{eqnarray*}
&&\frac{(-1)^{i-1}}{\prod_{1\le h<k}(\lambda_{m-h}-\lambda_{m-k})}
\\
&&\quad{}\times\sum_{j=0}^{i-(k+1)}
\frac{1}{(\lambda_{m-k}-\lambda
_{m-k-1-j})\prod_{0\le h\le i-(k+1),h\ne j}(\lambda_{m-k-1-j}-\lambda
_{m-k-1-h})}.
\end{eqnarray*}
By taking moduli and applying Lemma~\ref{lemma CG09} to the sum above,
we obtain
\[
\frac{(-1)^{i-1}}{\prod_{1\le h\le i,h\ne k}(\lambda_{m-k}-\lambda_{m-h})}.
\]
The result now follows from \eqref{eq:I-i} and \eqref{eq:rescale-p-i},
and from the fact that in the $K$-dimensional case, the probability of
going from $\mm$ to $\mm-\ii$,
conditional on $|\ii|$,
 is
\[
p\bigl(i_{1},\ldots,i_{K}; \mm,|\ii|\bigr) =
\frac{
{\binom{m_{1}}{i_{1}}}
\cdots
{\binom{m_{K}}{i_{K}}}
}{
{\binom{|\mm|}{|\ii|}}
}. 
\]
\upqed\end{pf*}
\end{appendix}

\section*{Acknowledgements}

The second author is supported by the European Research Council (ERC)
through StG \red{``N-BNP''} 306406.
The authors would like to thank Valentine Genon-Catalot and Aleksandar
Mijatovi\'c for helpful discussions.


%

\printhistory

\end{document}